\def\Bbb#1{{\bf #1}}

\def\fnote#1{\footnote}
\def\blacksquare{\hbox{\vrule width 4pt height 4pt depth 0pt}}

\def\cwleftpar#1#2{\leftskip #1 \rightskip #2 plus 1fill}
\def\cwrightpar#1#2{\leftskip #1 plus 1fill \rightskip #2}
\def\cwcenterpar#1#2{\leftskip #1 plus 1fill \rightskip #2 plus 1fill}
\def\cwfullpar#1#2{\leftskip#1\rightskip#2}

\def\cwoutdent#1#2{\llap{\hbox to #1{#2 \hss}}\ignorespaces}
\def\cwparbegin#1#2#3#4#5{
	\ifcase #1 \cwleftpar{#2}{#3}
	\or \cwrightpar{#2}{#3}
	\or \cwcenterpar{#2}{#3}
	\else \cwfullpar{#2}{#3}\fi
	\ifcase #4 \baselineskip = 1.5\baselineskip
	\or \baselineskip = 2\baselineskip
	\or \baselineskip = 3\baselineskip
	\else \baselineskip = 1\baselineskip\fi
	\ifdim #5 > 0in \else \noindent \fi
	\noindent\ignorespaces}
\documentclass{article}
\begin{document}
\advance \vsize by -1\baselineskip
\def\makefootline{
{\vskip \baselineskip \noindent \folio                                  \par
}}

\vspace*{2ex}

\noindent {\Huge Transports along Paths in\\[0.4ex] Fibre Bundles}\\[1.5ex]
\noindent {\Large II. Ties with the Theory of Connections and\\[0.5ex]
				Parallel Transports}

\vspace*{2ex}

\noindent Bozhidar Zakhariev Iliev
\fnote{0}{\noindent $^{\hbox{}}$Permanent address:
Laboratory of Mathematical Modeling in Physics,
Institute for Nuclear Research and \mbox{Nuclear} Energy,
Bulgarian Academy of Sciences,
Boul.\ Tzarigradsko chauss\'ee~72, 1784 Sofia, Bulgaria\\
\indent E-mail address: bozho@inrne.bas.bg\\
\indent URL: http://theo.inrne.bas.bg/$\sim$bozho/}

\vspace*{2ex}

{\bf \noindent Published: Communication JINR, E5-94-16, Dubna, 1994}\\[1ex]
\hphantom{\bf Published: }
http://www.arXiv.org e-Print archive No.~math.DG/0503006\\[2ex]

2003 PACS numbers: 02.40.Ma, 04.90.+e\\[2ex]

\noindent
{\small
The \LaTeXe\ source file of this paper was produced by converting a
ChiWriter 3.16 source file into
ChiWriter 4.0 file and then converting the latter file into a
\LaTeX\ 2.09 source file, which was manually edited for correcting numerous
errors and for improving the appearance of the text.  As a result of this
procedure, some errors in the text may exist.
}\\[2ex]

	\begin{abstract}
A review of the parallel transport (translation) in fibre bundles is
presented. The connections between transports along paths and parallel
transports in fibre bundles are examined. It is proved that the latter ones
are special cases of the former.
	\end{abstract}\vspace{3ex}

{\bf 1. INTRODUCTION}

\medskip
In the work [1], we have considered certain aspects of the general theory of transports along paths in arbitrary fibre bundles without investigating its ties with the ones of parallel transports and connections, which is the aim of the present paper.

Sect. 2 contains a review of the theory of parallel transports in fibre bundles adapted to suit our purposes. At first, are considered parallel transports generated by connections, after which attention is paid to the axiomatic approach to the concept of parallel transport. The main result of this paper, proved in Sect. 3, is that any parallel transport, axiomatically defined or generated by a connection, is a transport along paths satisfying certain additional conditions. Also some other ties between parallel transports and transports along paths are investigated. In Sect. 4 it is shown how linear transports along paths generated by derivations of tensor algebras [2] can be regarded as (axiomatically defined) parallel transports.

 In this paper, we shall use the following notation.

By $(E,\pi ,B)$ we denote an arbitrary fibre bundle with a base $B$, bundle space $E$ and projection $\pi :E  \to B [25,37,17]$. The fibres $\pi ^{-1}(x), x\in B$ are supposed to be homeomorphic.

An arbitrary real interval and a path in $B$ are denoted, respectively, by
$J$ and $\gamma :J  \to $B. If $B$ is a manifold, the tangent to $\gamma $
vector field is written as $dpt\gamma$. The path
$\bar{\gamma}\colon J  \to E$ is a lifting of $\gamma :J  \to B$ (resp.
through $u\in \pi ^{-1}(\gamma (J)))$ if $\pi \circ\bar{\gamma}=\gamma
($resp. and $(J)\ni u)$.

By $M, T_{x}(M)$ and $(T(M),\pi ,M)$ we denote, respectively, a
differentiable manifold, the tangent to it space at $x\in M$, and  the
tangent bundle to $M, T(M):=\cup_{x\in M}T_{x}(M)$.

Now for reference purposes, we shall summarize a certain material from [1].

A transport along paths I in $(E,\pi ,B)$ is a map $I:\gamma \to I^{\gamma
}$, where $I^{\gamma }:(s,t)\to I^{\gamma }_{s  \to t}, s,t\in J$ in which
the maps
\[
 I^{\gamma }_{s  \to t}:\pi ^{-1}(\gamma (s))  \to \pi ^{-1}(\gamma
(t)),\qquad (1.1)
\]
 satisfy the equalities ({\it id}$_{X}$is the identity map
of the set $X)$:
\[
 I^{\gamma }_{t  \to r}\circ I^{\gamma }_{s  \to t}=I^{\gamma }_{s  \to r},
\quad r,s,t\in J,\qquad (1.2)
 \]
 \[
  I^{\gamma }_{s  \to s}={\it id}_{\pi ^{-1_{(\gamma
 (s))}}},\quad   s\in J. \qquad (1.3)
 \]
  It is easily seen  that
\[
\Bigl( I^{\gamma }_{s\to t}\Bigr)^{-1}
=I^{\gamma }_{t  \to s}.\qquad (1.4)
\]
Important special classes of transports along paths are selected by one or
both of the conditions
\[
 I^{\gamma \mid J^\prime }_{s\to t}=I^{\gamma }_{s\to t}.
\quad s,t\in J^\prime ,\qquad (1.5)
\]
\[
 I^{\gamma \circ \tau }_{s\to t}=I^{\gamma }_{\tau (s)\to
\tau (t)},\quad  s,t\in J^{\prime\prime},\qquad (1.6)
\]
 where $\gamma \mid
J^\prime $ is the restriction of $\gamma $ on the subinterval $J^\prime $ of
$J$ and $\tau :J^{\prime\prime}  \to J$ is one-to-one map from the interval
$J^{\prime\prime}\subset {\Bbb R}$ onto J.

As for the transports along paths the types of the intervals $J, J^\prime $ and $J^{\prime\prime}$ are insignificant, in this work, for purposes which will be cleared up later, all real intervals are supposed to be closed, i.e. of type $[a,b]$ for some $a\le b, a,b\in {\Bbb R}$.

If the fibres $\pi ^{-1}(x), x\in B$ are differentiable manifolds (e.g. when
$(E,\pi ,B)$ is smooth [40]), one can consider the class of smooth transports
along paths obeying the condition
\[
 I^{\gamma }_{s\to t}\in Diff(\pi ^{-1}(\gamma (s)),\pi ^{-1}(\gamma
(t))), \quad  s,t\in J, \qquad (1.7)
\]
where Diff$(M,N)$ is the set of
diffeomorphisms from the manifold $M$ on the manifold N.

Any path $\gamma :[0,1]  \to B$ is called canonical (canonically defined). Its inverse path is $\gamma _{-}:=\gamma \circ \tau ^{c}_{-}$, where $\tau ^{c}_{-}:[0,1]  \to [0,1]$ is given by $\tau ^{c}_{-}(s):=1-s, s\in [0,1] [38,39]$. The (canonical) product of the paths $\gamma _{1},\gamma _{2}:[0,1]  \to B$ is the path $\gamma _{1}\gamma _{2}:[0,1]  \to B$ such that $(\gamma _{1}\gamma _{2})(s):=\gamma _{1}(2s)$ for $s\in [0,1/2]$ and $(\gamma _{1}\gamma _{2})(s)=\gamma _{2}(2s-1)$ for $s\in [1/2,1] [38,39]$.

\medskip
\medskip
 {\bf 2. REVIEW OF PARALLEL TRANSPORTS IN\\ FIBRE BUNDLES}

\medskip
This section contains a brief review of the concept "parallel transport" in fibre bundles. It will be a basis for comparison of the parallel transport with the transports along paths studied in $[1-3]$.

A common feature of most of the works $[4,5,6-8,12,15-18,21-35]$
dedicated to that problem is that in them as a basic object is taken the connection (in corresponding fibre bundles) and with its help the parallel transport is defined. In connection with this, one can distinguish the works $[9-11,13,14,19,20,36]$ in which as an initial (axiomatically given) object one takes the parallel transport which, in its turn, defines (and sometimes is identified with) the connection. It has to be noted that in these works in contrast to our considerations in $[1-3]$, main attention is paid to the dependence of the parallel transport on the curve (path) along which it is made.

\medskip
\medskip
 {\bf 2.1. PARALLEL TRANSPORT IN DIFFERENTIABLE\\
  FIBRE BUNDLES ENDOWED WITH CONNECTION}

\medskip
Let $(E,\pi ,B)$ be locally trivial differentiable and smooth (of class $C^{1})$ fibre bundle [17,25,37]. The fibre $\pi ^{-1}(\pi (u))$ through $u\in E$ is a manifold the tangent space of which at $u$ is denoted by $T^{v}_{u}(E):= :=T_{u}(\pi ^{-1}(\pi (u)))$. Evidently $T^{v}_{u}(E)\subset T_{u}(E)$. By definition $T^{v}_{u}(E)$ consists of vertical vectors [5,25,38,40].

{\bf Definition} ${\bf 2}{\bf .}{\bf 1} (cf. [5,35,40,41]). A$ connection
(of general form) in $(E,\pi ,B)$ is a smooth (of class $C^{1})
\dim(B)$-dimensional distribution $T^{h}(E):E\to T(E)$ such that the image
$T^{h}(E):u\to T^{h}_{u}(E)$ of $u\in E$ lies in $T_{u}(E)$ and is a direct
complement of $T^{v}_{u}(E)$ in $T_{u}(E)$, i.e.
\[
 T^{v}_{u}(E)\oplus T^{h}_{u}(E):=T_{u}(E),\qquad (2.1)
\]
 where $\oplus $ is
the direct sum sign. By definition $T^{h}_{u}(E)$ consists of horizontal
(with respect to the connection $T^{h}(E))$ vectors.

{\bf Definition 2.2} (cf. [5,21,35,40,41]). The smooth $(C^{1})$ path
$\bar{\gamma }:J\to E$ is horizontal (with respect to the connection
$T^{h}(E))$ if its tangent vector field $\dot{\bar{\gamma}}$ is horizontal,
i.e. if $\dot{\bar{\gamma}}(s)\in T^{h}_{\bar{\gamma}(s)}(E)$.

{\bf Remark.} In this definition and below we speak about smooth (or differentiable), of class $C^{1}$, paths as the corresponding generalizations for partially smooth paths are trivial.

{\bf Definition} ${\bf 2}{\bf .}{\bf 3} (cf. [5,21,35,40,41])$. The lift $\gamma \to \bar{\gamma }:J\to E$ of $\gamma :J\to B, ($resp. through $u\in E)$ is horizontal (with respect to the connection $T^{h}(E))$ if $\bar{\gamma }$ is a horizontal path (resp. through $u)$.

For defining the concept "parallel transport" in differentiable fibre bundles, of primary importance is, the question of the {\sl existence} of a {\sl unique} horizontal lift of a given path from the base
in the total space of the fibre bundle through any point above it. As has been pointed out in [40], p. 607, lemma $2, a$ sufficient condition for this is the fibre $\pi ^{-1}(b)$ for some $b\in B$ to be a {\sl compact} manifold. (Because of the local triviality, if this is so for some $b\in B$, the above property will be valid also for every $b\in $B.) The existence and uniqueness of the lifting mentioned are automatically fulfilled in principal fibre bundles $(G$-fibre bundles) [5,35,40,41] where they are assured from the additional requirement for the connection (called often $a G$-connection) to be invariant under the action of the structure group $G$ of the principal fibre bundle. (More strictly, if $R_{g}:E\to E$ is the right action generated from $g\in G$, then the connection $T^{h}(P)$ of the principal fibre bundle $(P,\pi ,B,G)$ is defined by the following three conditions: $1^{o}. T^{v}_{u}(P)\oplus T^{h}_{u}(P)=T_{u}(P), u\in P; 2^{o}. T^{h}_{u}(P)$ must depend differentiably on $u\in P; 3^{o}. R_{g*}T^{h}_{u}(P)=T^{h}_{R_{g}}(P)$, where $R_{g*}$is the differential of $R_{g}[5,10,40].)$ Another case for the existence of a unique horizontal lift of any path from the bases through every lying above it point is when the fibres of the fibre bundle are discrete (see $[17], pp. 75-76$ and [38], chapter III, lemma 15.1). The above pointed problem is considered from a general point of view in chapter III of the book [38] (see e.g. sections $12, 13, 15$, and 16 from it), where, in particular, are given the corresponding necessary and sufficient conditions for the existence of (maybe unique) lift of the pointed above form.

Let there be given a smooth fibre bundle $(E,\pi ,B)$ with connection $T^{h}(E)$ such that for every $C^{1}$path $\gamma :J\to B$ and every point $u\in \pi ^{-1}(\gamma (J))$ there exists a unique horizontal lift $\gamma \to \bar{\gamma }_{u}$of $\gamma $ through $u$, i.e. $\bar{\gamma }_{u}:J\to E, \pi \circ \bar{\gamma }_{u}=\gamma , u\in \bar{\gamma }_{u}(J)$ and   $_{u}(s)\in T^{h}_{\gamma _{u}}(E), s\in $J. (We will note that this assumption for the connection, without being mentioned, is unexplicitly  used in the considerations in sections 1 and 2 of [9].) Let $J=[a,b], a\le $b. Let us note that the considered below connections (and parallel transports) will be of Ehresmann's type (see [25], vol. 1, p. 314).

{\bf Definition} ${\bf 2}{\bf .}{\bf 4} (cf. [5,7,15-18,21-25,35,40-42])$.
The parallel transport (generated by $T^{h}(E))$ of the fibre $\pi
^{-1}(\gamma (a))$ onto the fibre $\pi ^{-1}(\gamma (b))$ along the path
$\gamma :[a,b]\to B$ is a diffeomorphism
\[
\varphi _{\gamma }:\pi ^{-1}(\gamma (a))\to \pi ^{-1}(\gamma (b)),\qquad
(2.2)
\]
such that if $u\in \pi ^{-1}(\gamma (a))$, then
$\varphi _{\gamma }:u\to \varphi _{\gamma }(u):=\bar{\gamma }_{u}(b)$, where
$\bar{\gamma }_{u}:[a,b]\to E$ is the unique horizontal lift of $\gamma $ in
$E$ through u.

{\bf Definition} ${\bf 2}{\bf .}{\bf 5} (cf. [5,9,40])$. The parallel transport defined by the connection $T^{h}(E)$ is a map $\varphi $ from the set of $C^{1}$paths in the base $B$ into the group Morf$(E,\pi ,B)$ of the bundle morphisms of $(E,\pi ,B)$, such that if $\gamma :[a,b]\to B$, then $\varphi :\gamma \to \varphi _{\gamma }\in  \in $Diff$(\pi ^{-1}(\gamma (a)),\pi ^{-1}(\gamma (b)))$, i.e. the image $\varphi _{\gamma }$is the defined by $T^{h}(E)$ parallel transport along $\gamma $, which is an element of the group Diff$(\pi ^{-1}(\gamma (a)),\pi ^{-1}(\gamma (b)))$ of diffeomorphisms between the fibres $\pi ^{-1}(\gamma (a))$ and $\pi ^{-1}(\gamma (b))$.

{\bf Proposition 2.1.} The parallel transport $\varphi $ has the following three basic properties:

a) Invariance under orientation preserving parameter changes, i.e. if
$\gamma :[a,b]\to B$ and $\tau :[c,d]\to [a,b], c\le d, a\le b$ is an
orientation preserving diffeomorphism, then
\[
\varphi _{\gamma \circ \tau }=\varphi _{\gamma }.\qquad (2.3)
\]

b) If $\gamma
_{-}:[0,1]\to B$ is the (canonical) inverse to $\gamma :[0,1]\to B$ path,
i.e. $\gamma _{-}(s)=\gamma (1-s), s\in [0,1]$, then
\[
 \varphi _{\gamma _{-}}=(\varphi _{\gamma })^{-1}.\qquad (2.4)
\]

c) If $\gamma _{1},\gamma _{2}:[0,1]\to B, \gamma _{1}(1)=\gamma _{2}(0)$
and $\gamma _{1}\gamma _{2}$is the (canonical) product of $\gamma _{1}$and
$\gamma _{2}($see Sect. 1), then
\[
 \varphi _{\gamma _{1}}=\varphi _{\gamma _{2}}\circ \varphi _{\gamma
_{1}}.\qquad (2.5)
\]

 {\bf Remark.} Because of (2.3) it is enough to consider
(2.4) and (2.5), as well as any other property of the parallel transport,
only for canonically defined paths in spite of the fact that they are valid
also for arbitrary ones.

{\bf Proof.} The proof of this proposition can be found, for example, in $[5,7-10,13,14,18,40].\blacksquare $

Here we shall drop the generality of the above considerations and till the end of the present section we will deal with the specific case of principal fibre bundles $[5,15,21-24]$.

At first, let us note that in principal fibre bundles the parallel transport $\varphi _{\gamma }$along $\gamma $ commutes with the right action $R_{g}, g\in G$ of the structure group $G$ on the total space of the fibre bundle [8,15], i.e. $\varphi _{\gamma }\circ R_{g}=R_{g}\circ \varphi _{\gamma }$for arbitrary path $\gamma $ and every $g\in $G.

On the other hand, in these fibre bundles the parallel transport can be
defined uniquely also by the right action of $G ($see e.g. [40], p. 632,
theorem 1 and [9]). In fact, let $\gamma :[a,b]\to B$ and $u\in \pi
^{-1}(\gamma (a))$. Then, due to the local triviality of $(E,\pi ,B) ($see
e.g. [43], p. 48), there exist a neighborhood $U$ of $\gamma (a)=\pi (u)$ and
a diffeomorphism $\psi :\pi ^{-1}(U)\to U  G, \psi (u):=(\pi (u),\chi (u))$,
where $\chi :\pi ^{-1}(U)\to G$ is right invariant, i.e.  $\chi (R_{g}u)=\chi
(u)g, g\in G$, and hence
\[
 (R_{g}\circ \psi ^{-1})(\pi (u),\chi (u))=R_{g}u =\psi ^{-1}(\pi
(R_{g}u),\chi (R_{g}u)) =\psi ^{-1}(\pi (u),\chi (u)g).
\]
 Denoting by  ${\bf e}$  the  unit  of  $G$,  we  find:
\[
\varphi _{\gamma }(u)
= (\varphi _{\gamma }\circ \psi ^{-1})(\gamma (a),\chi (u))
= (\varphi _{\gamma }\circ \psi ^{-1})(\gamma (a),{\bf e}\chi (u))
\]
\[
= (\varphi _{\gamma }\circ R_{\chi (u)}\circ \psi ^{-1})(\gamma
(a),{\bf e})
= (R_{\chi (u)}\circ \varphi _{\gamma }\circ \psi ^{-1})(\gamma (a),{\bf e})
=(R_{\chi (u)}\circ \psi ^{-1})(\gamma (b),  _{\gamma }),
\]
where  in  the  last equality  we  have used   the   fact   that $(\varphi
_{\gamma }\circ \psi ^{-1})(\gamma (a),{\bf e})\in \pi ^{-1}(\gamma (b))$ and
consequently there exists  a unique $ g _{\gamma }\in G$, which does not
depend on $u$ and is such that $(\varphi _{\gamma }\circ \psi ^{-1})(\gamma
(a),{\bf e})=\psi ^{-1}(\gamma (b), g_{\gamma })$. So, in principal fibre
bundles the  parallel transport $\varphi _{\gamma }$  along  $\gamma $  is
given   by   the equality $\varphi _{\gamma }(u)=(R_{\chi (u)}\circ \psi
^{-1})(\gamma (b), g_{\gamma })$. Hence, the definition of a parallel
transport $\varphi $ is equivalent to the definition of a map
$g\colon\gamma\mapsto g_\gamma$ from the set of $C^{1}$ paths in $B$ onto $G$
such that $ :\gamma \to g_{\gamma }$.

 {\bf Proposition 2.2.} The map $ g :\gamma \to   g_{\gamma }$has the
properties:
\[
g_{\gamma \circ \tau } = g_{\gamma },\qquad (2.6)
\]
\[
g_{\gamma_{-}}=(g_{\gamma })^{-1},\qquad (2.7)
\]
\[
g_{\gamma_{1}\gamma_2}
= g_{\gamma _{1}} g_{\gamma _{2}},\qquad (2.8)
\]
where $\gamma , \tau, \gamma_{-}, \gamma_{1}, \gamma_{2}$, and
$\gamma_{1}\gamma _{2}$ are defined in proposition 2.1.

{\bf Proof.} The equalities $(2.6)-(2.8)$ follow from the definition of $
_{\gamma }$and, respectively, the equalities $(2.3)-(2.5).\blacksquare $

Let us note that in some works, e.g. in [11,17,19,26,36,40], the third
property of a parallel transport is expressed not through the equality (2.5),
but by
\[
\varphi _{\gamma _{1}\gamma_2}=\varphi _{\gamma _{1}}\circ \varphi _{\gamma
_{2}},\qquad (2.5^\prime )
\]
 which, generally, is not true when using the
accepted by us notions: as $\varphi _{\gamma }$acts on the {\sl left}, then
(2.5) is valid but not (2.5  ). For (2.5  ) to be valid, as pointed out in
[17], p. 76, one has to change the orientations of $\gamma _{1}, \gamma
_{2}$, and $\gamma _{1}\gamma _{2}$; in fact, from (2.5) and (2.4) it follows
that
\[
\varphi _{(\gamma _{1}\gamma_2)_-}
=(\varphi _{\gamma _{1}})^{-1}
=(\varphi _{\gamma _{2}}\circ \varphi _{\gamma _{1}})^{-1}
=(\varphi _{\gamma _{1}})^{-1}\circ (\varphi _{\gamma _{2}})^{-1}
=\varphi _{(\gamma _{1})_-}\circ \varphi _{(\gamma _{2})_-},
\]
i.e.
\[
 \varphi _{(\gamma _{1}\gamma_2)_-}
=\varphi _{(\gamma _{1})_-}\circ \varphi _{(\gamma _{2})_-}.\qquad (2.9)
\]
 So, if we make the change $\varphi _{\gamma }\to \varphi
_{\gamma _{-}}=(\varphi _{\gamma })^{-1}, (2.5^\prime )$ will be valid but
not (2.5). Such is the case, for instance, in the works $[10, 11,19, 20]$ in
which the parallel transport is defined as the map $\varphi ^\prime :\gamma
\to \varphi $  :$=\varphi _{\gamma _{-}}:\pi ^{-1}(\gamma (b))\to \pi
^{-1}(\gamma (a))$, for which, due to $(2.9), (2.5^\prime )$ is true.

As regards the property (2.8) (in principal fibre bundles) in its right hand side the terms are written in a needed order as $  _{\gamma }$acts on the {\sl right} but not from the left as $\varphi _{\gamma }$.

At the end of this section, we shall stress the fact that the properties $(2.3)-(2.5)$ of the parallel transport $\varphi $ express their dependence on the curve of transport. From this viewpoint, there naturally arises the question of the "continuity" or "differentiability" (the "smoothness") of that dependence. The author knows two approaches to that problem. First, in the set of smooth (of class $C^{1})$ paths a topology is introduced (see e.g. [14] and [38], p. 104) which, in particular, may be generated by some metric (for a case of closed paths see [13]), which is used to study the smoothness of the map $\varphi :\gamma \to \varphi _{\gamma }$. And second, a (generally multidimensional) smooth deformation of $\gamma $ is made and the dependence of $\varphi _{\gamma }$on that deformation [9] is investigated, i.e. the class of homotopic with $\gamma $ paths connecting $\gamma (a)$ and $\gamma (b)$ is considered and the dependence of a parallel transport along these paths on the parameters of the used homotopy is investigated [38].

\medskip
\medskip
 {\bf 2.2. AXIOMATIC APPROACH TO THE PARALLEL\\
 TRANSPORT IN LOCALLY TRIVIAL FIBRE BUNDLES}

\medskip
The axiomatic definition of a parallel transport in locally trivial fibre bundles is based on the idea of a (diffeomorphic) mapping of the fibres of a given fibre bundle one onto another. More precisely, in the known to the author literature $[9-11,19,20,36,39,40]$ in which this question is set, it is put in the following way. Let $(E,\pi ,B)$ be a locally trivial fibre bundle and $x_{1},x_{2}\in $B. To any path $\gamma :J\to B$, where $J=[a,b]$, in the base $B$ connecting $x_{1}$and $x_{2}$, i.e. for which $\gamma (a)=x_{1}$and $\gamma (b)=x_{2}, a$ map (diffeomorphism$) \varphi _{\gamma }:\pi ^{-1}(x_{1})\to \pi ^{-1}(x_{2})$ is put into correspondence and the dependence of $\varphi _{\gamma }$on $\gamma $ is axiomatically defined. Namely, on $\varphi _{\gamma }$are imposed two kinds of restrictions. Firstly, these are conditions of a functional type defining the "change" of $\varphi _{\gamma }$when with the path $\gamma $ some operation is made (e.g. changing its orientation or its representation as a product of other paths). Secondly, in an appropriate way the "smoothness" of the map $\gamma \to \varphi _{\gamma }($conditions for smoothness) is defined. We shall note that the defined in this way parallel transport is sometimes called a global or an integral connection in the fibre bundle [9,19].

A scheme for solving the stated above problem for an axiomatic definition of the parallel transport in locally trivial fibre bundles has been introduced, maybe for the first time, in the work [19], after which, with little changes (following the context or using some features in different special cases (e.g. in principal or homogeneous (associated) fibre bundles)), it is repeated in other publications of the  same author [10,11,36].

\par
The above question, but in the "infinitesimal" case (the points $x_{1}$and $x_{2}$are infinitely near in a coordinate sense), is investigated in the works of G.F. Laptev (see [20] and the given therein references of the printed works of G.F. Laptev).

Ref. [9] contains a more general consideration of the problem, which is analogous to the one of Subsect. 2.1, but in $[9] a$ more general concept for connection ("infinitesimal nonlinear" connection) is used which is due to the replacement of the tangent spaces to the corresponding manifolds with the Grassmanian manifolds consisting of their one dimensional (linear) subspaces.

In [40], part II, sect. 24 the above question is described but, in fact, only a construction of a parallel transport by the method described in Subsect. 2.1 is made.

In the above sense, the defined in [39], sect. 3.2 transport along paths in an assembly of groups (a (flat) topological fibre bundle, the fibres of which are groups) is also a parallel transport.

Form here till the end of the present subsection we shall make comments on the axiomatic definition of the parallel transport in the mentioned above references and, in connection with our purposes attention will be paid mainly to the conditions of functional character.

Before going on, let us note that in the cited literature instead of an arbitrary closed interval $J=[a,b]$ the unit interval
${\bf I}=[0,1]$ is used, i.e. ${\bf I}=J\mid _{a=0,b=1}$. This is not important because of the invariance of the parallel transport under orientation preserving changes of the parameter of the paths along which it acts (see below $eq. (2.11))$.

Let $(E,\pi ,B)$ be a locally trivial smooth fibre bundle, $J=[a,b],
x_{1},x_{2}\in B, \gamma :J\to B$ be $a C^{1}$path, $\gamma (a)=x_{1}$and
$\gamma (b)=x_{2}$. The parallel transport in $(E,\pi ,B)$ is a map $\varphi
$ from the set of $C^{1}$paths in the base $B [38]$ onto the group
Morf$(E,\pi ,B)$ of bundle morphisms of $(E,\pi ,B) [17,22,25]$, such that
\[
\varphi :\gamma \to \varphi _{\gamma }\in
Diff(\pi ^{-1}(\gamma (a)),\pi ^{-1}(\gamma (b))).\qquad (2.10)
\]
 The first group of restrictions imposed on
$\varphi _{\gamma }$usually, contains $(2.3)-(2.5)$, i.e. it is wanted that
\[
\varphi _{\gamma \circ \tau }=\varphi _{\gamma },\qquad (2.11)
\]
\[
\varphi _{\gamma _{-}}=(\varphi _{\gamma })^{-1},\qquad (2.12)
\]
\[
\varphi _{\gamma _{1}}=\varphi _{\gamma _{2}}\circ \varphi _{\gamma
_{1}},\qquad (2.13)
\]
 where $\gamma , \tau , \gamma _{1}, \gamma _{2}$ and
$\gamma _{1}\gamma _{2}$ are defined in proposition 2.1.

The conditions $(2.11)-(2.13)$, which generally are independent, are
postulated, for example, in [10,11,19,36], where instead of $\varphi _{\gamma
}\sigma ^{\gamma }:=\Gamma (\gamma ):=\varphi _{\gamma _{-}}$is used, as a
consequence of which (2.13) is written in the form (2.5  ) (with $\sigma
^{\gamma }$instead of $\varphi _{\gamma }$, and, besides, the paths $\gamma $
and $\gamma ^\prime :=\gamma \circ \tau $ are called equivalent, which is
denoted by $\gamma\sim \gamma^\prime $, and (2.11) is written as $\sigma
^{\gamma }=\sigma ^{\gamma ^\prime }$ for $\gamma\sim \gamma ^\prime )$.

In [40] the restrictions $(2.11)-(2.13)$ are mentioned but, in fact, they are not used for an axiomatic construction of parallel transports.

In [9], attention is paid uniquely to the condition (2.13) which taken together with the corresponding condition for smoothness defines therein $\varphi $ as an integral connection of the fibre bundle. As in this work the full proofs of the stated there propositions are not given, part of which are not correct (e.g. the existence of a unique lift is supposed (see Subsect. 2.1); something which generally is not true (see e.g. [38,40])), the author of the present text was not able to re-establish them to an end, so it is not clear whether $(2.11), (2.12)$ or some other restrictions on $\varphi $
are used unexplicitly in [26].

Usually, as a consequence of other restrictions (resp. independently) (see
e.g. $[11,19,20]) \varphi $ satisfies (resp. on $\varphi $ is imposed) the
restriction
\[
\varphi _{\gamma _{a}} ={\it id}_{\pi ^{-1_{(x}}}, \gamma _{a}:\{a\}\to
\{x_{a}\}, x_{a}\in B, a\in {\Bbb R},\qquad (2.14)
\]
 i.e. to the degenerated
into a point path there  corresponds (resp. to correspond) the identity map
of the fibre over that point.

For example, in $[11,19] \gamma _{t}=\gamma \mid [0,t], t\in [0,1]$ is put
to be the restriction of $\gamma $ on $[0,t]\subset [0,1]$ and it is required
that   $i$  $_{0}\varphi _{\gamma _{t}}=${\it id}$_{\pi ^{-1_{(\gamma
(0))}}}(a$ functional condition) and that the principal part of the deviation
of $\varphi _{\gamma _{t}}$from {\it id}$_{\pi ^{-1_{(\gamma (0))}}}$should
depend smoothly on $\gamma _{t}$and   $_{t}($condition for smoothness) from
where, evidently, follows (2.14). On the contrary, if (2.14) is taken as a
base, then the first of these restrictions will be a consequence from the
condition for smoothness (which, in fact, needs a concrete and strict
formulation $(cf. [9]))$.

{\bf Definition 2.6.} The map $\varphi :\gamma \to \varphi _{\gamma }$, where $\varphi _{\gamma }$satisfies $(2.10)-(2.14)$, is called an axiomatically defined parallel transport.

 {\bf Remark.} In [40], p.~608, $\varphi$ is called an abstract connection.

From the described here approach to the parallel transport a little aside are the investigations of G.F. Laptev (see [20] and the references in it) due to their coordinate and local (or strictly - infinitesimal) character. As a consequence of this, the functional conditions and the conditions for smoothness (differentiability) of a parallel transport are given in a unified way (see [20], p. $46-47)$, not sharply separately as in our text or in [9]. From the above conditions in [20], p. 47 (see therein condition $c))$ only (2.14) is given, but $(2.11)-(2.13)$ therein are a consequence from the explicit coordinate and infinitesimal form of a parallel transport. Besides, in $[20] \varphi _{\gamma _{-}}$is used instead of $\varphi _{\gamma }$.

As has already been said above, the second group of restrictions imposed on the map (2.1) are the conditions for smoothness. They are defined $\varphi :\gamma \to \varphi _{\gamma }$as a continuous or differentiable (from some class $C^{k}, k=1,\ldots  ,\infty ,\omega )$ function of $\gamma $.

In the approach used in [20] these conditions are reduced to the requirement for analyticity of the principal linear part of an
explicit coordinate expression for the transport from the final point of a transport (see [20], p. 47, condition $d))$.

In the works of U.G. Lumiste [10,11,19,36] the question of smoothness of $\varphi :\gamma \to \varphi _{\gamma }$is, in fact, replaced with the requirement for continuous differentiability (smoothness) of the map $t\to \varphi _{\gamma _{t}}, t\in [0,1], \gamma _{t}:=\gamma \mid [0,t]$ with, maybe, some modifications depending on the concrete case under consideration, as is, for example, in [36], p.206, condition $\sigma 3$ where the concrete properties of the homogeneous fibre bundles are used. This condition for smoothness may be put in the first of the types described at the end of Subsect. 2.1 as it uses the topology of the real line (instead of the one in the set of smooth paths in $B)$.

We shall especially mention the work [9] where the important role of the conditions for smoothness is stressed and they themselves, in the considered there cases, are formulated strictly and clearly.

At the end of this section we shall only mention that there also exist a third group of conditions which sometimes are imposed on the map (2.1) and which are connected with the concrete structure of the investigated fibre bundles. They usually define the "intercommunications" of the map (2.1) with the (structural) group of transformations acting in the fibre bundle. Typical examples of this are the conditions $\sigma 2$ and $\sigma 4$ from [36], p. $205-206$ which concern homogeneous fibre bundles and the condition $\varphi _{\gamma }\circ R_{g}=R_{g}\circ \varphi _{\gamma }$for commutation of $\varphi _{\gamma }$with the right action $R_{g}, g\in G$ of the structure group $G$ in the case of principal fibre bundles [8,15].

\medskip
\medskip
 {\bf 3. THE AXIOMATICALLY DEFINED PARALLEL TRANSPORT\\
 AS A SPECIAL CASE OF TRANSPORTS ALONG PATHS}

\medskip
Before comparing a parallel transport with transport along paths we have to note the following. The axiomatically defined parallel transport is considered usually, along canonically given paths $\gamma :[0,1]\to B$, which is significant when defining explicitly the canonically inverse path $\gamma _{-}$and the canonical product of two paths (see Sect. 1 and [1], Sect. 3). Because of the invariance under parameter changes of the parallel transport (see Subsect. 2.2), this restriction is not essential and it is a question of convenience and easiness in the corresponding investigations. This
circumstance shows that the parallel transport must be compared not with the general transport along arbitrary paths, but with transports $I^{\gamma }$along the $\gamma :J\to B$, where $J$ is a closed interval, i.e. $J=[a,b]$. The importance of this restriction comes from the fact that, in the general case, the transports along paths are not invariant under parameter changes, i.e. they do not satisfy (1.6), so they can explicitly depend on the path of transport.

Let I be a transport along paths in the fibre bundle $(E,\pi ,B)$ and
$\gamma :[a,b]\to $B. To I we assign a map $\varphi :\gamma \to \varphi
_{\gamma }$, defined by
\[
 \varphi _{\gamma }:=I^{\gamma }_{a\to b}:\pi ^{-1}(\gamma (a))\to \pi
^{-1}(\gamma (b)).\qquad (3.1)
\]

 {\bf Lemma 3.1.} If $I^{\gamma }$ is a
transport along $\gamma $ satisfying additional conditions (1.5) and (1.6),
then the map $\varphi :\gamma   \to \varphi _{\gamma }$ defined by (3.1)
satisfies the equalities $(2.11)-(2.13)$ and
\[
 I^{\gamma }_{s\to t}=\varphi _{\gamma o\tau ^{J}_{t}}\circ (\varphi _{\gamma
o\tau ^{J}_{s}})^{-1}, s,t\in J=[a,b],\qquad (3.2)
\]
 where $\tau ^{J}_{s}:[a,b]\to [a,s], s\in [a,b]$ are for $s>a$ arbitrary
orientation preserving diffeomorphisms depending on $\gamma $ through the
interval J.

{\bf Proof.}  Firstly, we shall prove equality (3.2). Using sequentially
(1.2), (1.4), (1.5), (1.6) and (3.1), we get:
\[
 I^{\gamma }_{s\to t}
=I^{\gamma }_{a\to t}\circ I^{\gamma }_{s\to a}
=I^{\gamma }_{a\to t}\circ {\bigl(}I^{\gamma }_{a\to s}{\bigr)}^{-1}
=I^{\gamma \mid [a,t]}_{a\to t}\circ
\Bigl(I^{\gamma \mid [a,s]}_{a\to s}\Bigr)^{-1}
\]
\[
 =I^{\gamma \mid [a,t]}_{\tau ^{J_{(a)\to \tau }}_{t}}\circ
\Bigl(I^{\gamma \mid [a,s]}_{\tau ^{J_{(a)\to \tau }}_{s}} \Bigr)^{-1}
=I^{\gamma \circ \tau ^{J}_{t}}_{a\to b}\circ
\Bigl( I^{\gamma \circ \tau ^{J}_{s}}_{a\to b  } \Bigr)^{-1}
=\varphi _{\gamma \circ \tau ^{J}_{t}}\circ
 (\varphi _{\gamma \circ\tau ^{J}_{s}})^{-1}.
\]
 The property (2.11) follows from the equality (1.6): if $\tau :[c,d]\to
[a,b]$ is an orientation preserving diffeomorphism, which, in particular,
means $\tau (c)=a$ and $\tau (d)=b$, then from (1.6) and (3.1), we get
$\varphi _{\gamma \circ \tau }=I^{\gamma \circ \tau }_{c\to d}=I^{\gamma
}_{\tau (c)\to \tau (d)}=I^{\gamma }_{a\to b}=\varphi _{\gamma }$.

The property (2.12) is a consequence of (1.6) in which, because of $\gamma
_{-}:=\gamma \circ \tau ^{c}_{-}, \gamma :[0,1]  \to B$, we have to put $\tau
=\tau ^{c}_{-}($see Sect. 1). Under these assumptions, from $(3.1), (1.6)$
and (1.4) we get
\[
 \varphi _{\gamma _{-}}=I^{\gamma _{-}}_{0\to 1}=I^{\gamma }_{\tau
^{c_{(0)\to \tau }}_{-}}=I^{\gamma }_{1\to 0}=(\varphi _{\gamma })^{-1}.
\]

 The property (2.13) is a consequence of [1], proposition 3.4 (see therein
$eq. (3.4))$ in the case of a canonical choice of a parameter $\chi $, i.e.
(see Sect. 1 and [1]) for $\chi =\chi ^{c}:=(0,1,1/2;\tau ^{c}_{1},\tau
^{c}_{2})$ with $\tau ^{c}_{1}:s\to 2s, s\in [0,1/2]$ and $\tau ^{c}_{2}:s\to
2s-1, s\in [1/2,1]. ($It should be noted that the proof of proposition 3.4 of
[1] essentially uses the condition (1.5).) Then, from (3.1) and $eq. (3.4)$
from [1], we get
\[
 \varphi _{\gamma _{1}}=I^{\gamma _{1}}_{0\to 1}=I^{\gamma _{2}}_{0\to \tau
_{2}}\circ I^{\gamma _{1}}_{\tau _{1}}=I^{\gamma _{2}}_{0\to 1}\circ
I^{\gamma _{1}}_{0\to 1}=\varphi _{\gamma _{2}}\circ \varphi _{\gamma _{1}},
\]
 with $\gamma _{1}\gamma _{2}:={\bigl(}\gamma _{1}\gamma _{2}
\bigr)_{\chi ^{c}}.\blacksquare $

{\bf Lemma 3.2.} The defined by (3.1) map $\varphi :\gamma \to \varphi
_{\gamma }$for any transport along paths I has the property (2.14).

{\bf Proof.} If $\gamma _{a}:\{a\}\to \{x_{a}\}, x_{a}\in B, a\in {\Bbb R}$,
then from (3.1) and (1.3), we get $\varphi _{\gamma _{a}}=I^{\gamma
_{a}}_{a\to a}=${\it id}$_{\pi ^{-1_{(\gamma }}}=${\it id}$_{\pi
^{-1_{(x}}}.\blacksquare $

  {\bf Theorem 3.1.} If I is a smooth transport along paths, i.e.
 \[
I^{\gamma }_{a  \to b}\in Diff(\pi ^{-1}(\gamma (a)),\pi ^{-1}(\gamma(b))),
\quad \gamma :J  \to B,\ J=[a,b],
 \qquad (3.3)
\]
 having the properties (1.5) and
(1.6), then the defined by (3.1) map $\varphi :\gamma \to \varphi _{\gamma
}$is an axiomatically defined parallel transport. Vice versa, if $\varphi $
is an axiomatically defined parallel transport, then the map (3.2), in which
$\tau ^{J}_{s}:J  \to [a,s], s\in J$ are arbitrary orientation preserving
diffeomorphisms and $\gamma :J  \to B$, defines a smooth transport along
paths $I:\gamma \to I^{\gamma }, I^{\gamma }:(s,t)\to I^{\gamma }_{s  \to
t}$satisfying the additional conditions (1.5) and (1.6).

{\bf Remark.} If $\varphi $ is an axiomatically defined parallel transport,
then, because of the properties of $\tau ^{J}_{s}, s\in J$, we can replace in
$(2.11) \gamma $ with $\gamma \mid [a,s]$ and put in it $\tau =\tau
^{J}_{s}$. In this way, we obtain $\varphi _{\gamma \mid [a,s]}=\varphi
_{(\gamma \mid [a,s])\circ \tau ^{J}_{s}}=\varphi _{\gamma \circ \tau
^{J}_{s}}$as $(\gamma \mid [a,s])\circ \tau ^{J}_{s}=\gamma \circ \tau
^{J}_{s}$. Therefore, (3.2) is now equivalent to
\[
 I^{\gamma }_{s  \to t}=\varphi _{\gamma \mid [a,t]}\circ (\varphi _{\gamma
\mid [a,s]})^{-1}, s,t\in J=[a,b].\qquad (3.4)
\]

 {\bf Proof.} The first part
of the theorem is a consequence of lemmas 3.1 and 3.2, definition 2.6 and the
fact that now (2.10) is, due to (3.1), another form of (3.3).

On the contrary, let $\varphi $ be an axiomatically defined parallel transport (see definition 2.6).

If in theorem 3.1 of [1], we put $Q=\pi ^{-1}(\gamma (a))$ and $F^{\gamma }_{s}= =$  $\varphi _{\gamma o\tau ^{J}_{s}}$  $^{-1}:\pi ^{-1}(\gamma (s))\to \pi ^{-1}(\gamma (a)) (\tau ^{J}_{s}(a)=a, \tau ^{J}_{s}(b)=s)$, we see that the map (3.2) is a transport along $\gamma $ from $s$ to t. So, $I:\gamma \to I^{\gamma }$, where $I^{\gamma }:(s,t)\to I^{\gamma }_{s  \to t}$, is a transport along paths.

 The smoothness condition (3.3) follows from (2.10) and (3.2).

To prove the equalities (1.5) and (1.6) for the transport
along paths I, we shall use the following lemma which will be proved below after this proof.

{\bf Lemma 3.3.} If $\varphi $ is an axiomatically defined parallel
transport, then the maps (3.4) (or equivalently (3.2)) admit the
representation
\[
I^{\gamma }_{s  \to t}=(\varphi _{\gamma \mid
[\min(s,t),\max(s,t)]})^{\epsilon (s,t)}
=
\Big\{	\begin{array}{ll}
\varphi_{\gamma \mid [s,t]}&  \ \ for\ s\le t \\
\varphi_{\gamma \mid [t,s]}& \ \  for\ {s\ge t},
	\end{array}
\qquad (3.5)
\]
 where $\epsilon (s,t):=+1$ for $s\le t$ and $\epsilon (s,t):=-1$ for $s>t
($or $s\ge t)$.

From (3.5), because of $(\gamma \mid J^\prime )\mid J^\prime =\gamma \mid
J^\prime $ for any subinterval $J^\prime \subseteq J$, it immediately follows
\[
  I^{\gamma }_{s  \to t}=I^{\gamma \mid [\min(s,t),\max(s,t)]}_{s  \to
t}\qquad (3.6)
\]
 which by [1], proposition 2.3 is equivalent to (1.5).

If $\tau :J^{\prime\prime}  \to J$ is an orientation preserving
diffeomorphism, then $(\gamma \circ \tau )\mid [r,s]=(\gamma \mid [\tau
(r),\tau (s)]\circ \tau $ for every $r,s\in J^{\prime\prime}$ such that $r\le
$s. Combining this equality with (3.5), letting $s,t\in J^{\prime\prime},
\lambda :=\min(s,t)$ and $\mu := :=\max(s,t)$, and using (2.11), we get:
\[
 I^{\gamma \circ \tau }_{s  \to t}=(\varphi _{(\gamma \circ \tau )\mid
[\lambda ,\mu ]})^{\epsilon (s,t)}=(\varphi _{(\gamma \mid [\tau (\lambda
),\tau (\mu )])\circ \tau })^{\epsilon (s,t)}
\]
\[
 =(\varphi _{(\gamma \mid [\tau (\lambda ),\tau (\mu )])})^{\epsilon
(s,t)}=(\varphi _{(\gamma \mid [\tau (\lambda ),\tau (\mu )])})^{\epsilon
(\tau (s),\tau (t))}=I{ } ^{\gamma }_{\tau (s)  \to \tau (t)}
\]
 as $r\le s$ leads to $\tau (r)\le \tau (s), r,s\in J^{\prime\prime}.\blacksquare $

 The proof of lemma 3.3 is based on

{\bf Lemma 3.4.} If $\varphi $ is an axiomatically defined parallel
transport and $\gamma :J  \to B$, then
\[
 \varphi _{\gamma \mid [s,t]}\circ
\varphi _{\gamma \mid [r,s]}
=\varphi _{\gamma \mid [r,t]}, \qquad  for\ r\le s\le t,\ r,s,t\in J.
 \qquad (3.7)
\]

 {\bf Proof.} Let $\tau _{1}:[0,1]  \to [r,s]$ and
$\tau _{2}:[0,1]  \to [s,t]$ be orientation preserving diffeomorphisms.
Evidently, also such is the map $\tau :[0,1]  \to [r,t]$, defined by $\tau
(\lambda ):=\tau _{1}(2\lambda )$ for $\lambda \in [0,1/2]$ and $\tau
(\lambda ):=\tau _{2}(2\lambda -1)$ for $\lambda \in [1/2,1]$. Using (2.11),
the definition of the (canonical) product of paths (see Sect. 1), and (2.13),
we find:
\[
  \varphi _{\gamma \mid [s,t]}\circ \varphi _{\gamma \mid [r,s]}=\varphi
_{\gamma \circ \tau _{1}}\circ \varphi _{\gamma \circ \tau _{2}}=\varphi
_{(\gamma \circ \tau _{1}}=\varphi _{\gamma \circ \tau }=\varphi _{\gamma
\mid [r,t]}.\blacksquare
\]

{\bf Proof of lemma 3.3.} Combining (3.7) and (3.4) for $a\le s\le t\le b$,
we get
\[
 I^{\gamma }_{s  \to t}=\varphi _{\gamma \mid [a,t]}\circ (\varphi _{\gamma
\mid [a,s]})^{-1}=\varphi _{\gamma \mid [s,t]}\circ \varphi _{\gamma \mid
[a,s]}\circ (\varphi _{\gamma \mid [a,s]})^{-1}=\varphi { } _{\gamma \mid
[s,t]}
\]
and for $a\le t\le s\le b$, we obtain
\[
 I^{\gamma }_{s  \to t}=\varphi _{\gamma \mid [a,t]}\circ (\varphi _{\gamma
\mid [a,s]})^{-1}=\varphi _{\gamma \mid [a,t]}\circ (\varphi _{\gamma \mid
[t,s]}\circ \varphi _{\gamma \mid [a,t]})^{-1}=(\varphi _{\gamma \mid
[t,s]})^{-1} \blacksquare
\]

Theorem 3.1 is a strict expression of the statement  that  the
axiomatically defined parallel transport is a special case of tran  sports
along paths in fibre bundles, and that any  transport  along paths satisfying
certain additional conditions,  namely  (1.5)  and (1.6), defines an
axiomatically defined  parallel  transport.  This theorem also expresses a
one-to-one correspondence between  axioma  tically defined parallel
transports  and  transports  along  paths obeying the conditions (1.5) and
(1.6). Speaking  more  freely,  we can say that according to it a transport
along paths is an  axioma  tically defined parallel transport if and only if
it satisfies  the additional conditions (1.5) and (1.6).

{\bf Proposition 3.1.} If a transport along paths I (resp. axiomatically defined parallel transport $\varphi )$ defines through (3.1) (resp. (3.2)) the axiomatically defined parallel transport $\varphi  ($resp. transport along paths I), then the generated by $\varphi  ($resp. I) by means of (3.2) (resp. (3.1)) transport along paths (resp. axiomatically defined parallel transport) coincides with the initial transport along paths I (resp. the axiomatically defined parallel transport $\varphi )$.

{\bf Proof.} Let $^\prime I ($resp. $^\prime \varphi )$ be the generated by
$\varphi  ($resp. I) transport along paths (resp. axiomatically defined
parallel transport). Using (3.1) and (3.2), we find
\[
 ^\prime I^{\gamma }_{s  \to t}
=\varphi _{\gamma \circ \tau ^{J}_{t}}\circ
\varphi ^{-1}_{\gamma \circ \tau ^{J}_{s}}
=I^{\gamma \circ \tau ^{J}_{t}}_{a  \to b}\circ
\Bigl( I^{\gamma \circ \tau ^{J}_{s}}_{a  \to b  } \Bigr)^{-1}
\]
\[
=I^{\gamma }_{\tau^J_t(a) \to \tau^J_t(b) }  \circ
 I^{\gamma }_{\tau^J_s(b) \to \tau^J_s(a) }
=I^{\gamma }_{a\to t}\circ I^{\gamma }_{s  \to a}=I^{\gamma }_{s  \to t}
\]
(resp. $^\prime \varphi _{\gamma }=I^{\gamma }_{a  \to b}=\varphi _{\gamma
\circ \tau ^{J}_{b}}\circ \varphi ^{-1}_{\gamma \circ \tau ^{J}_{a}}=\varphi
_{\gamma }\circ ${\it id}$_{\pi ^{-1_{(\gamma (a))}}}=\varphi _{\gamma
}).\blacksquare $

\medskip
\medskip
{\bf 4. THE GENERATED BY DERIVATIONS OF\\ TENSOR ALGEBRAS
 TRANSPORTS ALONG PATHS AS\\ PARALLEL TRANSPORTS
 IN TENSOR BUNDLES}

\medskip
In this section, by $\eta $ we denote a $C^{1}$ path in the manifold $M$ such
that $\eta :[a,b]\to M$ for a definite $a\le b, a,b\in {\Bbb R}$.

Let $S$ be an $S$-transport along paths (in the tensor algebra over $M) [2]$.

{\bf Definition 4.1.} The $S$-parallel transport associated with the
$S$-transport $S$ is a map $\varphi $ from the set of $C^{1}$paths in $M$
into the set of bundle morphisms of the tensor bundles over these paths such
that
\[
 \varphi :\eta \to \varphi _{\eta }:=S^{\eta }_{a\to b}:T_{\eta (a)}(M)  \to
T_{\eta (b)}(M), \qquad  a\le b,\qquad (4.1)
\]
 where $T_{x}(M)$ is the tensor algebra
at $x\in $M. The map $\varphi _{\eta }$will be called an $S$-parallel
transport along (the path$) \eta $.

{\bf Lemma 4.1.} If $\varphi $ is the $S$-parallel transport generated by an
$S$-transport $S, \gamma :J\to M$ and $s,t\in J$, then
\[
 S^{\gamma }_{s\to t}
=
\Big\{
	\begin{array}{ll}
 {\varphi_{\eta }}&\ \ \eta=\gamma|[s,t] \quad for\ s\le t \\
 (\varphi _{\eta })^{-1}& \ \ \eta=\gamma|[t,s] \quad for\ t\le  s
	\end{array}
. \qquad (4.2)
\]

 {\bf Proof.} (4.2) follows from (4.1) and (1.4), as any
$S$-transport has this property (see $[2], eq. (2.10)$ and also [3], Sect.
2).\blacksquare

Between the $S$-transports and $S$-parallel transports there exists one important difference. Namely, the $S$-transport along $\gamma :J\to M$ does not use the natural order of the real numbers which defines a definite orientation on the interval $J$, while in the definition (4.1) of an $S$-parallel transport this order is used explicitly $(a\le b)$. The last fact is the reason for the appearance of two different cases $(s\le t$ and $s\ge t)$ in (4.2). This fact also reflects the difference between (1.6) (or (4.4)) and (4.5) (see below proposition 4.2).

  {\bf Proposition 4.1.} If $\eta _{a}:\{a\}\to \{m_{a}\}, a\in {\Bbb R}$ and
$m_{a}\in M$, then
\[
 \varphi_{\eta_a} ={\it id}_{\pi^{-1}(m_a)} .\qquad (4.3)
\]

 {\bf Proof.} (4.3) follows directly from (4.1) for $b=a$ and (1.3) (see also
[2], definition 2.1).\blacksquare

{\bf Proposition 4.2.} Let $\eta :[a,b]\to M, \tau :[a^\prime ,b^\prime ]\to [a,b]$ be a diffeomorphism and the $S$-transport $S^{\eta }$along $\eta $ be invariant under
the change $\tau $ of the parameterization of $\eta $, i.e. $(cf. (1.6))$
\[
 S^{\gamma \circ \tau }_{s\to t}=S^{\gamma }_{\tau (s)\to \tau (t)},
\quad s,t\in [a^\prime ,b^\prime ].\qquad (4.4)
\]
 Then, for the $S$-parallel transport
$\varphi $, corresponding to $S$, there holds
\[
 \varphi _{\eta \circ \tau }=\varphi _{\eta },\quad for\ \tau (a^\prime
)=a,\qquad (4.5a)
\]
\[
 \varphi _{\eta \circ \tau }=(\varphi _{\eta })^{-1},  \qquad
for\ \tau (a^\prime )=b,\qquad (4.5b)
\]
 i.e. the $S$-parallel transport is
invariant under orientation preserving change of the parameterization (case
(4.5a)), but when this change does not preserve the orientation it is
replaced by its inverse map (case (4.5b)).

{\bf Proof.} Using successively $(4.1), (4.4)$ and (4.2) for $\gamma =\eta
$, we find
\[
 \varphi _{\eta \circ \tau }
=S^{\eta \circ \tau }_{a^\prime \to b^\prime }
=S^{\eta }_{\tau (a^\prime )\to \tau (b^\prime )}
=
\Big\{
	\begin{array}{ll}
(\varphi_\eta & \ \ for\ {\tau (a^\prime )\le \tau (b^\prime )} \\
(\varphi_\eta)^{-1} & \ \ for\  {\tau (a^\prime )\ge \tau (b^\prime )}  ,
	\end{array}
\]
 which, due to that $\tau $ is a diffeomorphism, is equivalent to
(4.5).\blacksquare

{\bf Proposition 4.3.} Let $\eta _{-}:=\eta \circ \tau ^{c}_{-}$be the
canonically inverse path to $\eta :[a,b]\to M$, i.e. $\tau ^{c}_{-}:[a,b]
\to [a,b], \tau ^{c}_{-}(s):=a+b-s, s\in [a,b] (cf$. Sect. 1 and [38,39]). If
(4.4) holds for changing the orientation maps $\tau $ for some $S$-transport,
then for the corresponding to it $S$-parallel transport $\varphi $ there
holds the equality
\[
 \varphi _{\eta _{-}}=(\varphi _{\eta })^{-1}.\qquad (4.6)
\]

 {\bf Proof.} This
result is a corollary from $(4.1), (4.4), (1.4)$ and the inequality $(\tau
^{c}_{-})^{-1}(a)\ge (\tau ^{c}_{-})^{-1}(b)$ as, by definition $\tau
^{c}_{-}$ changes the orientations. $ Eq. (4.6)$ also follows from (4.5b) for
$\tau =\tau ^{c}_{-}.\blacksquare $

{\bf Proposition 4.4.} Let $\eta _{1}\eta _{2}$be the (canonical) product of
the paths $\eta _{h}:[0,1]\to M, h=1,2, \eta _{1}(1)=\eta _{2}(0) ($see Sect.
1). If an $S$-transport defining the $S$-parallel transport $\varphi $
satisfies (4.4) for preserving the orientations $\tau $ and (1.5), then
\[
 \varphi _{\eta _{1}}=\varphi _{\eta _{2}}\circ \varphi _{\eta _{1}}.\qquad
(4.7)
\]

 {\bf Proof.} Putting $\tau _{1}(s)=2s, s\in [0,1/2]$ and $\tau
_{2}(s)=2s-1, s\in [1/2,1]$ and using sequentially $(4.1), (1.2), (1.5)$ and
(4.4), we get:
\[
 \varphi _{\eta _{1}}= S^{\eta _{1}}_{0\to 1} = S^{\eta _{1}}_{1/2}\circ
S^{\eta _{1}}_{0\to 1/2}= S^{(\eta _{1}}_{1/2\to 1}
\]
\[
  \circ  S^{(\eta _{1}}_{0\to 1/2}= S^{\eta _{2}}_{1/2\to 1}\circ  S^{\eta
_{1}}_{0}= S^{\eta _{2}}_{\tau _{2}}
\]
\[
  \circ  S^{\eta _{1}}_{\tau _{1}}= S^{\eta _{2}}_{0\to 1}\circ  S^{\eta
_{1}}_{0\to 1}= \varphi _{\eta _{2}}\circ  \varphi _{\eta _{1}}.\blacksquare
\]

In propositions $4.2, 4.3$ and 4.4 one essentially uses the acceptance for the validity of (4.4). This is not random as the equality (4.4) expresses the invariance  (under certain conditions) under the changes of parameterization of an $S$-transport's path, and all (parallel) transports (see Sect. 2) known to the author and used in the mathematical and physical literature possess this property.

From the above-said it is clear that under sufficiently general and "reasonable" conditions an $S$-parallel transport satisfies all basic (functional) conditions characterizing the parallel transport when it is axiomatically  described (see Sect. 2.2). Namely, this is the reason for calling the map (4.1) an $S$-parallel transport: it is a "parallel transport" acting in the tensor spaces over a differentiable manifold and it is generated by derivation of the tensor algebra over the manifold. More precisely, from the above results and definition 2.6, we derive

{\bf Proposition 4.5.} The $S$-parallel transport generated by an $S$-transport along paths satisfying along them (1.5) and (1.6) is the axiomatically defined parallel transport.

The next proposition expresses some properties of the $S$-parallel transports which are specific of them as "parallel transports" in tensor bundles.

{\bf Proposition 4.6.} Any $S$-parallel transport $\varphi _{\eta }$along a path $\eta :[a,b]\to M$ possesses the properties:

a) Linearity: if $\lambda ^\prime ,\lambda ^{\prime\prime}\in {\Bbb R}$ and
$T^\prime $and $T^{\prime\prime}$ are tensors at $\eta (a)$, then:
\[
 \varphi _{\eta }(\lambda ^\prime T^\prime +\lambda
^{\prime\prime}T^{\prime\prime})=(\lambda ^\prime )\varphi _{\eta }(T^\prime
+(\lambda ^{\prime\prime})\varphi _{\eta }(T^{\prime\prime});\qquad (4.8)
\]

 b) Term by term action on tensor products: if A and $B$ are arbitrary
tensors at $\eta (a)$, then
\[
\varphi _{\eta }(A\otimes B)=(\varphi _{\eta }(A))\otimes (\varphi _{\eta
}(B));\qquad (4.9)
\]

 c) Commutativity with the contraction operator $C$:
\[
 \varphi _{\eta }\circ C-C\circ \varphi _{\eta }=0;\qquad (4.10)
\]

 d) An identical action on scalars: if $\lambda \in {\Bbb R}$, then
\[
 \varphi _{\eta }(\lambda )=\lambda .\qquad (4.11)
\]

{\bf Proof.} Equalities $(4.8)-(4.11)$ follow directly from definition 4.1
and, respectively, the properties $(2.3)-(2.5)$ and (2.13) of [2] of the
$S$-transports.\blacksquare

\medskip
\medskip
 {\bf 5. CONCLUSION}

\medskip
 The main result of this work is that the theory of transports along paths in
fibre bundles is sufficiently general and includes as its special case the
theory of parallel transports, and also, consequently, the connection theory.
An essential role, as we saw, in comparing these theories was played by the
additional condition (1.6). The transports along paths satisfying it depend
in fact not on the {\it path}  of transport $\gamma :J  \to B$ but on the
{\it curve} of transport, i.e. on the whole class of paths $\{\gamma \circ
\tau \}$ in which $\tau $ is $a 1:1$ map of ${\Bbb R}$-intervals onto J.
Because of the practical importance of (1.6), we shall consider it below in
the most used case, the one of linear transports in vector bundles [3].

Let $L$ be a linear transport in the vector bundle $(E,\pi ,B) [3]$. If
$\tau :J^{\prime\prime}  \to J$ is 1:1 map, then $eq. (1.6)$ reads
\[
 L^{\gamma o\tau }_{s\to t}=L^{\gamma }_{\tau (s)\to \tau (t)}, s,t\in
J^{\prime\prime}.\qquad (5.1)
\]

 Let a field of bases in $E$ be fixed along
$\gamma :J  \to B$ in which $H:(s,t;\gamma )\to H(s,t;\gamma ), s,t\in J$ and
$\Gamma _{\gamma }(s):=(\partial H(s,t;\gamma )/\partial t)\mid _{t=s}$be,
respectively, the matrix and the matrix of the coefficients of $L [3]$. Let
${\cal D}^{\gamma }$be the generated by $L$ derivation along $\gamma $ and
$({\cal D}^{\gamma }\sigma )(s)=:{\cal D}^{\gamma }_{s}\sigma $ for $a
C^{1}$section $\sigma $ of $(E,\pi ,B) ($see [3], eqs. (4.2) and (4.3)).

{\bf Proposition 5.1.} The condition (5.1) is equivalent to any of the
following three equalities:
\[
H(t,s;\gamma \circ \tau )=H(\tau (t),\tau (s);\gamma ), s,t\in
J^{\prime\prime},\qquad (5.2)
\]
\[
\Gamma _{\gamma \circ \tau }(s)
=  \frac{d\tau(s)}{ds} \cdot \Gamma _{\gamma }(\tau (s)),
 \qquad  s\in J^{\prime\prime},\qquad (5.3)
\]
\[
{\cal D}^{\gamma \circ \tau }_{s}
=  \frac{d\tau(s)}{ds} \cdot {\cal D}^{\gamma }_{\tau (s)},
 \qquad s\in J^{\prime\prime}.\qquad (5.4)
\]

 {\bf Proof.} The equivalence of (5.1) and
(5.2) is a corollary of the definition of $H ($and the linearity of $L$; see
[3], Sect. 2). Eqs. (5.3) and (5.4) are equivalent because of the connection
(4.7) or $(4.14^\prime )$ from [3] between $\Gamma _{\gamma }(s)$ and ${\cal
D}^{\gamma }_{s}$.

So, it remains to prove the equivalence between (5.2) and (5.3).

Differentiating (5.2) with respect to $s$ and using $\Gamma _{\gamma }(s):=
:=(\partial H(s,t;\gamma )/\partial t)\mid _{t=s}$, we get (5.3). On the
contrary, if (5.3) holds, then using the same equality, the representation
$H(t,s;\gamma )= =F^{-1}(t;\gamma )F(s;\gamma )$ for some matrix function $F
($see [3], proposition 2.4) and $dF^{-1}/ds=-F^{-1}(dF/ds)F^{-1}$, we easily
obtain:
\[
\frac{d}{ds}
[H(\tau (t),\tau (s);\gamma )H^{-1}(t,s;\gamma \circ \tau )]
\]
\[
 = H(\tau (t),\tau (s);\gamma )
\Bigl[ \frac{d\tau(s)}{ds} \cdot \Gamma _{\gamma }(\tau (s)) - \Gamma
_{\gamma \circ \tau }(s)
\Bigr]  H^{-1}(t,s;\gamma \circ \tau ) = 0.
\]
  From this, due to $H(s,s;\gamma )={\Bbb I} ($see $[3], eq. (2.12))$ and
$(t=s  \Leftrightarrow  \tau (t)=\tau (s))$, we derive
(5.2).\blacksquare

If $B$ is a manifold, evident examples of linear transports along paths satisfying (5.3), and hence (5.1), are the ones characterized by the coefficients given by $[3], eq. (5.1)$ and, in particular, the parallel transports generated by linear connections.

The definition of a parallel transport in principal or associated fibre
bundles by the map $\gamma \to g_{\gamma }\in G$ (see Subsect. 2.1) is
widely used in the physical literature devoted to gauge theories $[13,14,
29-32,43-45]$. In them, the parallel transport is given globally  through an
ordered (called also $P-, T$-, or chronological) exponent [14,40,45] along
$\gamma $, i.e.
$\gamma \to g_{\gamma }= Pexp \int A_{i}dx^{i}$, where
$A_{i}$ are the components of the connection form (or, in physical language,
the gauge potentials). So, locally along a path $\gamma $ connecting the
infinitesimally near points $x$ and $x+dx$ it is defined by the expansion
$g_{\gamma }\mid _{x,x+dx}={\Bbb I}+A_{i}dx^{i}[11,43]$.

If $\gamma $ is a closed path (a contour) passing through $x\in B ($in the
physical literature such a path is called a loop), then the quantity $W(\gamma ,x):=$Pexp  $A_{i}dx^{i}$is called a Wilson loop $[29-32]$ and in accordance with the above considerations it uniquely defines the parallel transport from $\pi ^{-1}(x)$ onto $\pi ^{-1}(x)$, i.e. of the fibre over $x$ onto itself. The importance of Wilson's loops is in that their set $\{W(\gamma ,x): \gamma :[a,b]\to B, \gamma (a)=\gamma (b), x\in \gamma ([a,b])\}$, which is a nonabelian group and is a representation of the group of loops, contain all the information for the considered gauge theory $[13,29-32,43-45]$.

\medskip
\medskip
 {\bf ACKNOWLEDGEMENT}

\medskip
This research was partially supported by the Fund for Scientific Research of Bulgaria under contract Grant No. $F 103$.

\medskip
\medskip
 {\bf REFERENCES}

\medskip
1.  Iliev B.Z., Transports along paths in fibre bundles. General theory, JINR Communication $E5-93-299$, Dubna, 1993.\par
2.  Iliev B.Z., Parallel transports in tensor spaces generated by derivations of tensor algebras, JINR Communication JINR $E5-93-1$, Dubna, 1993. \par
3.  Iliev B.Z., Linear transports along paths in vector bundles. I. General theory, JINR Communication JINR $E5-93-239$, Dubna, 1993.\par
4.  Hawking S.W., G.F.R. Ellis, The large scale structure of space-time, Cambridge Univ. Press, Cambridge, 1973.\par
5.  Kobayashi S., K. Nomizu, Foundations of differential geometry, vol.1, Interscience publishers, New-York-London, 1963.\par
6.  Sachs R.K., $H{\bf .} Wu$, General Relativity for Mathematicians, Springer-Verlag, New York-Heidelberg-Berlin, 1977.\par
7.  Nash C., S.Sen, Topology and Geometry for physicists, Academic Press, London-New York, 1983.\par
8.  Nomizu K., Lie groups and differential geometry, The mathematical Society of Japan, 1956.\par
9.  Nikolov P., On the correspondence between infinitesimal and integral description of connections, ICTP, Internal Report, $IC/81/196$.\par
10.  Lumiste U.G., Connection theory in fibre bundles, In: Science
review, sec. Mathematics: Algebra. Topology. Geometry. 1969, VINITI, Moscow, $1971, 123-168 ($in Russian).\par
11.  Lumiste U.G., Parallel transport; Connections on manifolds; Connection object; Connection form; Connections in fibre bundles, Articles in: Mathematical Encyclopedia, vol.4, Soviet encyclopedia, Moscow, 1984 (in Russian).\par
12.  Warner F.W., Foundations of differentiable manifolds and Lie groups, Springer Verlag, New York-Berlin-Heidelberg-Tokyo, 1983.\par
13.  Durhuus B., J.M. Leinaas, On the loop space formulation of gauge theories, CERN, preprint TH 3110, Geneva, 1981.\par
14.  Khudaverdian O.M., A.S. Schwarz, A new comments on the string representation of gauge fields, Phys. Lett. ${\bf B}, {\bf 9}{\bf 1}$, No.$1, 1980, 107-110$.\par
15.  Bishop R.L., R.J. Crittenden, Geometry of Manifolds, Academic Press, New York-London, 1964.\par
16.  Yano K., M. Kon, Structures on Manifolds, Series in Pure Mathematics, vol. 3, World Scientific Publ. Co., Singapore, 1984.\par
17.  Steenrod N., The topology of fibre bundles, 9-th ed., Princeton Univ. Press, Princeton, $1974 ( 1-st ed. 1951)$.\par
18.  Sulanke R., P. Wintgen, Differentialgeometrie und Faserb\"undel, VEB Deutscher Verlag der Wissenschaften, Berlin, 1972.\par
19.  Lumiste U.G., To the foundations of global connection theory, In: Scientific writings of Tartu state university, No.150, Works on mathematics and mechanics IV, Tartu, $1964, 69-108 ($in Russian).\par
20.  Ostianu N.M. V.V. Rizhov, P.I. Shveikin, Article on the scientific investigations of German Fedorovich Laptev, In: Works of Geometrical seminar, vol.4, VINITI, Moscow, $1973, 7-70 ($in Russian).\par
21.  Choquet-Bruhat Y. et al., Analysis, manifolds and physics, North-Holland Publ.Co., Amsterdam, 1982.\par
22.  Husemoller D., Fibre bundles, McGrow-Hill Book Co., New York-St. Louis-San Francisco-Toronto-London-Sydney, 1966.\par
23.  Mishtenko A.S., Vector fibre bundles and their applications, Nauka, Moscow, 1984 (in Russian).\par
24.  Hermann R., Vector bundles in mathematical physics, vol.I, Benjamin W.A., Inc., New York, 1970.\par
25.  Greub W., S. Halperin, R. Vanstone, Connections, Curvature, and Cohomology, vol.1, vol.2, Academic Press, New York and
London, $1972, 1973$.\par
 26.  Atiyah M.F., $K$-theory, Harvard Univ., Cambridge, Mass., 1965. \par
27.  Karoubi M., $K$-theory. An Introduction, Springer-Verlag, Berlin-Heidelberg-New York, 1978.\par
28.  Dandoloff R., W.J. Zakrzewski, Parallel transport along a space curve and related phases, J. Phys. {\bf A}: Math. Gen., ${\bf 2}{\bf 2}{\bf ,} 1989, L461-L466$.\par
29.  Azam M., Gauge-invariant objects from Wilson loops, Phys. Rev. ${\bf D}, {\bf 3}{\bf 5}$, No.$6, 1987, 2043-2046$.\par
30.  Jevicki A., B. Sakita, Collective approach to the large $N$ limit: Euclidean field theories, Nucl. Phys. ${\bf B}{\bf ,} {\bf 1}{\bf 8}{\bf 5}, 1981, 89-100$.\par
31.  Polyakov A.M., Gauge fields as rings of glue, Nucl. Phys. ${\bf B}{\bf ,} {\bf 1}{\bf 6}{\bf 4}, 1979, 171-188$.\par
32.  Makeenko Yu.M., A.A. Migdal, Quantum chromodynamics as dynamics of loops, Nucl. Phys. ${\bf B}{\bf ,} {\bf 1}{\bf 8}{\bf 8}, 1981, 269-316$.\par
33.  Tamura I, Topology of foliations, Mir, Moscow, 1979 (in Russian; translation from Japanese).\par
34.  Hicks N.J., Notes on Differential Geometry, D. Van Nostrand Comp., Inc., Princeton, 1965.\par
35.  Sternberg S., Lectures on differential geometry, Chelsea Publ. Co., New York, 1983.\par
36.  Lumiste U.G., Homogeneous fibre bundles with connection and their immersions, In: Works of the Geometrical seminar, vol.1, VINITI, Moscow, $1966, 191-237 ($in Russian).\par
37.  Hartman Ph., Ordinary differential equations, John Wiley \& Sons, New-York-London-Sydney, 1964.\par
$38. Hu$ Sze-Tsen, Homotopy Theory, Academic Press, New York-London, 1959.\par
39.  Viro O.Ya., D.B. Fuks, I. Introduction to homotopy theory, In: Reviews of science and technology, sec. Modern problems in mathematics. Fundamental directions, vol.24, Topology-2, VINITI, Moscow, $1988, 6-121 ($in Russian).\par
40.  Dubrovin B.A., S.P. Novikov, A.T. Fomenko, Modern geometry, Nauka, Moscow, 1979 (in Russian).\par
41.  Lichnerowicz A., Th\^eorie globale des connexions et des gropes $D  $Holonomie, Roma, Edizione Crem\'onese, 1955.\par
42.  Mangiarotti L., M. Modugno, Fibred Spaces, Jet Spaces and Connections for Field Theories, In: Proceedings of the International meeting on Geometry and Physics, Florence, October $12-15.1982, ed$. M. Modugno, Pitagora Editrice, Bologna, $1983, 135-165$.\par
43.  Yang C.N., Integral formalism for gauge fields, Phys. Rev. Lett., {\bf 33}, No.$7, 1974, 445-447$.\par
44.  Gambini R., A. Trias, Geometrical origin of gauge theories, Phys. Rev. ${\bf D}{\bf ,} {\bf 2}{\bf 3}$, No.$2, 1981, 553-555$.\par
45.  Aref$^{,}$eva I.Ya., On the integral formulation of gauge theories, preprint No.480, Wroclaw Univ., Wroclaw, 1979.

\newpage
\noindent
 Iliev B. Z.\\[5ex]

\noindent
 Transports along Paths in Fibre Bundles\\
 II. Ties with the Theory of Connections and Parallel Transports\\[5ex]

\medskip
A review of the parallel transport (translation) in fibre bundles is
presented. The connections between transports along paths and parallel
transports in fibre bundles are examined. It is proved that the latter ones
are special cases of the former.\\[5ex]

\medskip
\medskip
\medskip
\medskip
 The investigation has been performed at the Laboratory of Theoretical
Physics, JINR.

\end{document}